\documentclass{article}

\begin{document}

\bigskip
\bigskip

Tsemo Aristide

College Boreal

951, Carlaw Avenue,

M4K 3M2

 Toronto, ON Canada.

tsemo58@yahoo.ca

\bigskip
\bigskip

{\bf Scheme theory for groups  groups and Lie algebras.}

\bigskip

{\bf Introduction.}

\medskip

Let ${\cal L}= {\cal F}\bigcup {\cal P}\bigcup {\cal C}$ be a
first order language, ${\cal F}$ is the set of symbols of
operations of arities $n_F$, ${\cal P}$ is the set of symbols of
predicates, and ${\cal C}$ a set of constants. Examples of such languages are:
 the language of groups, the operations are the
multiplication, the inversion, and the neutral is a constant.
  The language of Lie algebras, the operations here
are the addition, the substraction, the Lie bracket, constants
$0,1$.

An ${\cal L}$ structure is defined by a set $M$, and for each
symbol $F$ in ${\cal F}$, a function $F^M:M^{n_F}\rightarrow M$,
where $n_F$ is the arity of $F$, a constant $c_M$ for each element
$c$ in ${\cal C}$; if ${\cal L}$ is the language of groups, a structure is
a group. If ${\cal L}$ is the language of Lie algebras, a
structure is a Lie algebra.

Let $X=\{x_1,..,x_n\}$ a set of variables, the free ${\cal L}$
structure $T_{\cal L}[X]$ is a ${\cal L}$ structure universal in
the sense that for any ${\cal L}$-structure $M$,  any morphism of
sets $X\rightarrow M$, extends uniquely to a morphism of ${\cal
L}$-structures $T_{\cal L}(X)\rightarrow M$. Example, if ${\cal
L}$ is the language of groups, $T_{\cal L}(X)$ is the free group
generated by $X$.
The formulas or equation on $T_{\cal L}(X)$, are defined
recursively with elements of $T_{\cal L}(X)$ and logic operations.
If ${\cal L}$ is the language of groups, we can always write a set
of formulas $(S=1), S\subset T_{\cal L}(X)$, if ${\cal L}$ is the
language of Lie algebras, we can write $(S=0)$.

In the sequel, we suppose that the language used are algebraic,
these are the language which define groups, Lie algebras,
algebras,...

Let $M$ be an ${\cal L}$-structure. A $M$-solution of an equation
$(S=1)$ is a morphism $f:T_{\cal L}(X)\rightarrow M$ such that
$f(S)=1$. The set of solutions of the equation $(f=1)$ is denoted
$V(S)$. These sets generate the closed subsets of a topology
called the Zariski topology.

The algebraic geometry on groups or (Lie) algebras is  the study
of the geometry of closed subsets for the Zariski topology above.
The  motivation is to provide a geometric interpretation of subsets of groups
and Lie algebras defined by algebraic equations. Like the commutator of a finite subset of a group.
The main results of this theory have been stated
for groups by Baumslag, Myasnikov, and Remeslennikov, and for Lie
algebras by Daniyarova, Kazachkov, Remeslennikov.

The previous authors, focus on  the coordinates systems, the purpose of this paper
is to develop a free coordinates approach inspired by the scheme theory defined by Grothendieck.

\bigskip

{\bf Algebraic geometry of groups.}

\medskip

Firstly, we are going to study the algebraic geometry of groups. Let $C$ be the category
of groups, and $G$ an object of $C$. We denote by $G\mid C$ the comma category whose objects
are morphisms $\phi_H:G\rightarrow H$, such an object will usually be denoted $(H,\phi_H)$.
Recall that a morphism $f:(H,\phi_H)\rightarrow (H',\phi_{H'})$ is a morphism of groups
$f:H\rightarrow H'$ such that $f\circ \phi_{H}=\phi_{H'}$ We will be mainly concerned in the
full subcategory $C(G)$ of $G\mid C$ whose objects are the objects $(H,\phi_H)$ such
that $\phi_H$ is injective.

\medskip

{\bf Definition.}

Let $x$ be an element of the object $(H,\phi_H)$ of $C(G)$, we denote by $G(x)$, the subgroup of $H$
generated by $\{gxg^{-1}, g\in G\}$.

The element $x$ is invertible if and only if $G(x)\cap \phi_H(G)\neq 1$.

A non trivial element $x$ of $H$ is a divisor of zero if and only if
 there exists a non trivial element $y$ of $H$  such that the group of commutators
  $[G(x),G(y)]$ of elements of $G(x)$ and $G(y)$ is trivial.

A $G$-domain is an object $(H,\phi_H)$ of $C(G)$ which does not have zero divisors.

An normal subgroup $P$ of $H$ is a prime ideal if and only if:

- $P\cap \phi_H(G)=1$

- $H/P$ a $G$-domain. This is equivalent to saying that for every element $x,y\in H$,
  such that  $[G(x),G(y)]\subset P$, $x\in P$ or $y\in P$.

\medskip

Let $I$ be a normal subgroup of $H$, we denote by $V(I)$ the set of prime normal subgroups
which contains $I$. We denote by $Spec(H)$ the set of prime normal subgroups of $H$. We call it
an affine $G$-scheme.

\medskip

{\bf Remarks.}

If $I\cap G$ is different of $\{1\}$, then $V(I)$ is empty.

If $x$ is a divisor of zero, for each $G$-automorphism $h$, $h(x)$ is a divisor of zero.

\medskip

{\bf Proposition.}

{\it Let $I$ and $J$ be two normal subgroups of $H$, $V([I,J])=V(I)\bigcup V(J)$

For any family of normal subgroups $(I_a)_{a\in A}$, let $I$ be the subgroup of $H$
generated by the groups $(I_a)_{a\in A}, V(I)=\bigcup_{a\in A}V(I_a)$.}

\medskip

{\bf Proof.}

Let $P$ be a prime subgroup which is an element of $V(I)\bigcup V(J)$, $P$ contains $I$
or $J$, this implies that $P$ contains $[I,J]$. Thus $V(I)\bigcup V(J)\subset V([I,J])$.

Let $P$ be an element of $V([I,J])$, suppose that there exists  elements $x\in I, y\in J$.
such that $x,y$ are not in $P$. The group $G(x)\subset I$, and $G(y)\subset J$,
thus the subgroup $[G(x),G(y)]$ of $H$ is contained in $[I,J]$. This implies that $x$ is in $P$
or $y$ is in $P$.

Let $I$ be the subgroup of $H$ generated by the family $(I_a)_{a\in A}$. Consider a prime $P$ which
contains $I$. This implies that $P$ contains $I_a$ for every element $a\in A$. We
deduce that $V(I)\subset \bigcap_{a\in A}V(I_a)$. Let $P$ be a prime in $\bigcap_{a\in A}V(I_a)$,
$P$ contains $I_a, a \in A$, this implies that $P$ contains $I$.

\medskip

The previous proposition shows that there exists a topology on $Spec(H)$ for which the closed
subsets can be written $V(I)$, where $I$ is a normal subgroup of $H$.

Remark that if $G$ is commutative, $Spec(H)$ is always empty. Since for every element $x$
of $G$, $[G(x),G(x)]=1$.

\bigskip

{\bf Localization and structural sheaves.}

\medskip

Let $(H,\phi_H)$ be an object of $C(G)$, we are going to define three sheaves on $Spec(H)$:

\medskip

Let $U$ be an open subset of $Spec(H)$, we define $L_{Spec(G)}(U)$ to be the set of applications
$f:U\rightarrow \prod_{P\in U} H/P$,
such that for every $P\in U$, there exists an open subset $V$ containing $P$, an element
$f_V\in H$ such that for every $Q\in V$, $f(Q)$ is the image of $f_V$ by the quotient map $u_Q:H\rightarrow H/Q$.

\medskip

Let $H$ be a group, $Z(H)$ is the algebra of $H$: it is the $Z$-ring generated by $\{1_h, h\in H\}$.
Its multiplication structure is defined by $1_h1_{h'}=1_{hh'}$.
We denote by $A_{Spec(H)}(U)$ the space of functions $f:U\rightarrow \prod_{P\in U}Z(H/P)$
such that for every $P\in U$, there exists a neighborhood $V$ of $P$, an element $f_V\in Z(H)$
such that for every $Q\in V$, we have $f(Q)=p_Q(f_V)$, where $p_Q:Z(H)\rightarrow Z(H/Q)$
is the canonical morphism.

\medskip

To define the last sheaf, we recall the following result on the localization of modules of associative algebras
due to Schonfield:

\medskip

{\bf Theorem See Lidia Angeleri Hugel, and Maria Archetti.}

{\it Let $R$ be an associative ring, and $\Sigma$ a set of maps between finitely generated
projective right modules, then there exists an associative $R_{\Sigma}$, and a morphism
$\alpha:R\rightarrow R_{\Sigma}$ which is $\Sigma$-inverting; this is equivalent to saying that for every map
$f:M\rightarrow M'\in \Sigma, f\otimes Id_{R_{\Sigma}}:M\otimes_R\otimes R_{\Sigma}
\rightarrow M'\otimes_RR_{\Sigma}$ is an isomorphism. The morphism $\alpha$ is universal in
the sense that if $h:R\rightarrow S$ is another $\Sigma$-invertible morphism, then
there exists a morphism $h':R_{\Sigma}\rightarrow S$, such that $h=h'\circ\alpha$.}

\medskip

Let $P$ be an element of $Spec(H)$, and $p_P:Z(H)\rightarrow Z(H/P)$ the natural
projection. For any subset $D$ of $Spec(H)$,
we denote by $\Sigma_D$ the subset of elements of $Z(H)$ such that for every
$f\in \Sigma_D$, and every elements $P$ of $D$, $p_P(f)\neq 0$.
 Each element of $\Sigma_D$ defines a morphism of
the right module $Z(H)$ by left multiplication. We denote by $Z(H)_D$ the localization
of $Z(H)$ by $\Sigma_D$. And by $\alpha_D$ the inverting morphism of $\Sigma_D$.

Let $U$ be an open subset of $Spec(H)$, we denote by $O_{Spec(H)}(U)$ the set of functions
$f:U\rightarrow \prod_{P\in U}Z(H)_P$, such that for every element $P\in U$, there exists
a neighborhood $V$ of $P$ in $U$, an element $f_V\in  Z(H)_V$ such that
for every element $Q$ in $V$, $f(Q)$ is the image of $f_V$ by the canonical morphism
$Z(H)_V\rightarrow Z(H)_P$ resulting from the universal property of $Z(H)_P$.

\bigskip

{\bf Proposition.}

{\it Let $(H,\phi_H)$ and $(H',\phi_{H'})$ be elements of $C(G)$, and $f:H\rightarrow H'$
 be a morphism, then $f$ induces a continuous map $f^*:Spec(H')\rightarrow Spec(H)$,
 a morphism of sheaves $f_L:L_{Spec(H)}\rightarrow L_{Spec(H')}$, morphisms
 of ringed spaces $f_O:(Spec(H),O_{Spec(H)})\rightarrow (Spec(H'),O_{Spec(H')})$,
$f_A:(Spec(H),A_{Spec(H)})\rightarrow (Spec(H'),A_{Spec(H')})$.}

\medskip

{\bf Proof.}
Let $P$ be an element of $Spec(H')$, we set $f^*(P)=f^{-1}(P)$. The morphism
$H/f^{-1}(P)\rightarrow H'/P$ is injective, it implies that $f^{-1}(P)$ is a point of $Spec(H)$,
since a sub $G$-group of a $G$-group which does not have zero divisors, does not
have zero divisors.

Let $I$ be a normal subgroup of $H'$, $f^{-1}(V(I))=V(f^{-1}(I))$. This implies that $f^*$
is continuous.

Let $U$ be an open subset of $Spec(H')$, we define
$f_O(U):O_{Spec(H)}({f^*}^{-1}(U))\rightarrow O_{Spec(H')}((U)$ as follows:
Remark that $f^{-1}(\Sigma_U)\subset \Sigma_{{f^*}^{-1}(U)}$. Thus $f$ induces
a morphism $Z(H)_{{f^*}^{-1}(U)}\rightarrow Z(H')_{U}$ which is $f_O(U)$.

The maps $f_L$ and $f_A$ are defined analogously.

\medskip

{\bf Definition.}

A group $G$-scheme is a ringed space $(X,O_X)$, such that every element of $X$ has a neighborhood
isomorphic to an affine $G$-scheme.

\bigskip

{\bf Representations and schemes.}

\medskip

Let $L$ be a $G$-domain, $\{1\}$ is a prime, it is the generic point of $Spec(H)$. Suppose that $H$ is a $G$-domain,
An $L$-point of $Spec(H)$ is a morphism $i_L:Spec(L)\rightarrow Spec(H)$. The morphism
$i_L$ is defined by a morphism of $G$-groups $h_P:H\rightarrow L$. We denote by $P=h_P^{-1}(1)$.
In particular, if $P$ is an element of $Spec(H)$ such that $H/P=L$, it is an $L$-point.

The purpose of this part is to show that any affine scheme can be realized as a set of representations:
There exists a group $L$, such that every element of $Spec(H)$ is an $L$-point. On this
purpose we show the following:

\medskip

{\bf Lemma.}

{\it Let $(L_i)_{i\in I}$ be a family of $G$-domains, the amalgamated sum $L=\coprod_GL_i$
is a $G$-domain.}

\medskip

{\bf proof.}

Since the morphism $\phi_i:G\rightarrow L_i$ is injective, the canonical imbedding
$i_{L_i}:L_i\rightarrow L$ is injective (See Serre Theorem 1 p.9). Let $x$ and $y$ be two elements of $L-G$. Suppose
that $[G(x),G(y)]=1$. If there exists an element $i\in I$ such that $x=i_{L_i}(x_i),
y=i_{L_i}(y_i)$ we deduce that $[G(x_i),G(y_i)]=1$ since $i_{L_i}$ is injective. We deduce that
$x=1$ or $y=1$ since $L_i$ is an $H$-domain.
If such an $i$ doesn't exist, necessarily $[G(x),G(y)]\neq 1$, since $L$ is the quotient
of the free product of $(L_i)_{i\in I}$ by the normal subgroup generated by
$i_{L_i}(g){i_{L_j}(g)}^{-1}, g\in G$.

\medskip

{\bf Theorem.}

{\it Let $(H,\phi_H)$ be an object of $C(G)$. There exists a group $L$, such that every point
of $Spec(H)$ is an $L$-point.}

\medskip

{\bf Proof.}

Let $L$ be the amalgamated sum of $H/P, P\in Spec(H)$. Then the previous lemma shows that $L$
is a $G$-domain. Let $l_P:H/P\rightarrow L$, be the canonical imbedding. The composition of the
canonical projection $H\rightarrow H/P$ with $l_P$ endows $P$ with the structure of an $L$-point.

\medskip

{\bf Schemes defined by finite groups.}

\medskip

We are going to study here the theory developed above for finite groups. We start by the following
result:

\medskip

{\bf Proposition.}

{\it Let $G$ ba finite group, there exists a group $L$ such that
for every $(H,\phi_H)$ of $C(G)$ where $H$ is finite, the points of $Spec(H)$ are $L$-points.}

\medskip

{\bf Proof.}

Let $n$ be an integer, the set of isomorphic classes of $G$-domain of cardinality inferior to $n$ is finite.
This implies that the set of isomorphic classes of finite $G$-domains is numerable. Let $(L_i)_{i\in I}$
be a set of $G$-domains such that every $G$-domain is isomorphic to a $L_i$. We can suppose
that $L_i$ is finite or numerable. We define the amalgamated sum $L=\coprod_GL_i{i\in I}$.
If $Spec(H)$ is a $G$-affine scheme, for every point $P$ of $Spec(H)$, $H/P$ is a finite $G$-domain.
There exists an imbedding $H/P\rightarrow L$ whose composition with $H\rightarrow H/P$
defines a $L$-point.

\medskip

Let $G$ be a group, and $Aut(G)$ the group of automorphisms of $G$. An inner automorphism
$i_g$ of $G$ is  defined by an element $g\in G$ such that for every element $x$ of $G$,
$i_g(x)=gxg^{-1}$. We call $In(G)$ the group of inner automorphisms $G$, and by $Out(G)$
the quotient $Aut(G)/Inn(G)$. For a large class of groups, $Out(G)=1$
in that situation, we say that $G$ is complete. For example if $G=S_n, n\neq 2,6$
$G$ is complete,..

\medskip

{\bf Proposition.}
{\it Suppose that $H$ is a finite $G$-domain where $Out(G)=1$ , then the normalizer of $G$
in $H$ is $G$.}

\medskip

{\bf Proof.}

Let $N(G)$ be the normalizer of $G$ in $H$. Consider the map: $N(G)\rightarrow Aut(G)$,
$n\rightarrow i(n)$, the restriction of the inner morphism defined by $n$ to $G$. If $N(G)$
is not equal to $G$, then the kernel of $i$ is not trivial. Let $x$ be a non trivial element
in the kernel of $i$, $x$ commutes with $G$. This implies that $[G(x),G(x)]=1$. Thus $x$ is a divisor of zero.
This is a contradiction.

\medskip

We give an example of domain:

\medskip

{\bf Proposition.}

{\it Let $G=S_n, n>4$ be the symmetric group, then $S_n$ is an $S_n$-domain; $S_{n+1}$
endowed with the $S_n$-structure defined by the canonical imbedding $S_n\rightarrow S_{n+1}$
is an $S_n$-domain.}

\medskip

{\bf Remark.}

Another topology can be defined on the category of groups as follows:

 Let $G$ be a group, we say that a normal subgroup $P$ of $G$ is prime ideal if and only if for every
 normal subgroups $I$, $J$ of $G$, $[I,J]\subset P$ if and only if $I\subset P$ or $J\subset P$.
 We denote by $Spec(G)$ the set of prime ideals of $G$.
 Let $I$ be a normal subgroup of $G$, we denote by $V(I)$ the set of prime ideals of $G$ which contain
 $I$. The sets $V(I)$, are the closed subsets of a topology on $Spec(H)$.

 The problem with this approach is due to the fact that it is not functorial: if $H$ is a finite
 commutative group, $Spec(H)$ is empty, if $G$ is a simple group, $Spec(G)$ contains
 only one element. They may exist an imbedding $H\rightarrow G$ which does not induce
 a morphism $Spec(G)\rightarrow Spec(H)$.

\medskip

{\bf The affine scheme associated with equations.}

\medskip

The motivation of the study of algebraic geometry for groups is the study of algebraic equations
in group theory. Let $G$ be a group, and $X=\{x_1,...,x_n\}$ a finite set. Consider $H=G*F(x_1,...,x_n)$ the
the free product of $G$ with the free group generated by $X$. Many interesting sets in group
theory can be expressed with ideals in $H$. For example, if $g$ is an element of $G$, we can study
the set $E(g)$ $n$-uples $(g_1,...,g_n)\in G^n$ such that $(g_ig=gg_i$.  This $n$-uples are defined by
the equations $gx_ig^{-1}{x_i}^{-1}=1, i-1,..,n$ which generates a normal subgroup $I(g)$ of $H$.
The elements of $E(g)$ can be identified with morphisms $H\rightarrow G$ whose kernel contains $I(g)$.
Thus there is a bijective correspondence between $E(g)$ and $H/I(g)$.
Such system of equations are studied by Baumslag and his coauthors.

To study the algebraic equations on $H$, we can also study the algebraic scheme defined by its $G$-structure.
prime elements, correspond to generic points.

\medskip

{\bf Algebraic geometry of Lie algebras.}

\medskip

We are going firstly to study,  the topology on the set of maximal ideals of a Lie algebra

Let $Lie$ be the category of Lie algebras and $Lie({\cal S})$  the comma category whose objects are morphisms
${\cal S}\rightarrow {\cal G}$. We study here the full subcategory $C({\cal S})$ of $C({\cal S})$
whose objects are injective morphisms $\phi_{\cal G}:{\cal S}\rightarrow {\cal G}$. We denote
this object by $({\cal G},\phi_{\cal G})$.
Recall that a morphism between $f:({\cal G},\phi_{\cal G})\rightarrow ({\cal H},\phi_{\cal H})$
is a morphism of Lie algebras $f:{\cal G}\rightarrow {\cal H}$ such that
$\phi_{\cal H}=f\circ \phi_{\cal G}$.

 For every element $x\in ({\cal G},\phi_{\cal G})$, we denote by ${\cal S}(x)$ the orbit of $x$ by ${\cal S}$.

The category $Lie({\cal S})$ has limits and colimits, since the category of Lie algebras has
limits and colimits. A direct construction of sum in $Lie({\cal S})$ can be done follows:
 Let $({\cal G},\phi_{\cal G})$ and $({\cal H},\phi_{\cal H})$ two objects $C({\cal S})$.
 The sum of $({\cal G},\phi_{\cal G})$ and
$({\cal H},\phi_{\cal H})$ is the pushout of ${\cal G}$ and ${\cal H}$ by $\phi_{\cal G}$ and
$\phi_{\cal H}$; It is the quotient of ${\cal G}\oplus{\cal H}$ by the ideal generated
by $\{\phi_{\cal G}(s)-\phi_{\cal H}(s), s\in {\cal S}\}$.

The product of $({\cal G},\phi_{\cal G})$ and $({\cal H},\phi_{\cal H})$ endowed with
 the diagonal action of ${\cal S}$ defines a product in $C({\cal S})$.

\medskip

{\bf Definition.}

A non zero element $x\in ({\cal G}-{\cal S})$ is a divisor of zero, if there exists a non zero element
$y\in {\cal G}-{\cal S}$ such that $[{\cal S}(x),{\cal S}(y)]=0$.

- A prime ideal of $({\cal G},\phi_{\cal G})$ is an ideal $P$ of ${\cal G}$ such that:

- $P\cap \phi_{\cal G}({\cal S})=\{0\}$

Let $\phi_P:{\cal G}\rightarrow {\cal G}/P$ be the canonical projection, $({\cal G}/P,\phi_P\circ\phi_{\cal G})$
does not have a divisor of zero.
This is equivalent to saying that for every elements $x,y$ of ${\cal G}$  such
that $[{\cal S}(x),{\cal S}(y)]\subset P$, $x\in P$ or $y\in J$.

\medskip

For each ideal $I$ of ${\cal G}$, we denote by $V(I)$ the set of prime
ideals which contain $I$, and by $Spec({\cal G})$ the set of prime ideals of ${\cal G}$.

\medskip

{\bf Proposition.}

{\it Let ${\cal G}$ be a Lie algebra, for every ideals $I,J$ of ${\cal G}$,
we have: $V([I,J])=V(I)\bigcup V(J)$.

For any family of ideals $(I_p)_{p\in P}$, we have:
 $V(\oplus_{a\in A}I_a)=\bigcap_{a\in A}V(I_a)$.}

\medskip

{\bf proof.}

We firstly show that $V([I,J])=V(I)\bigcup V(J)$.
Let $P\in V(I)\bigcup V(J)$.
Since $P$ contains $I$ or $J$, it contains $[I,J]$. This implies that
$V(I)\bigcup V(J)\subset V([I,J])$.

Let $P$ be an element of $V([I,J])$, Suppose that there exists $x\in I, y\in J$
which are not elements of $P$. Since $I$ and $J$ are ideals, $x\in I$, and $y\in J$,
we deduce that ${\cal S}(x)\in I, {\cal S}(y)\in J$, and $[{\cal S}(x),{\cal S}(y)]\subset P$.
we deduce that $x\in P$ or $y\in P$ since $P$ is a prime ideal.

Now we show that $V(\oplus_{a\in A}I_a)=\bigcap_{a\in A}V(I_a)$.
Let $P$ be an element of $V(\oplus_{a\in A}I_a)$, $P$ contains $\oplus_{a\in A}I_a$.
This implies that $I_a\subset P$, for every $a\in A$, it results that $P\in V(I_a)$.
Thus $P\in \bigcap_{a\in A}V(I_a)$.

Conversely, let $P$ be an element of $\bigcap_{a\in A}V(I_a)$, $I_a\subset P$
for every $a\in A$.
This implies that $\oplus_{a\in A}I_a\subset P$.

\bigskip

{\bf Localization and the structural sheaf.}

\medskip

 Let $({\cal G},\phi_{\cal G})$ be an element of $C({\cal S})$,
  we denote by $E({\cal G})$ the enveloping of $G$. Recall that that we can imbed
  ${\cal G}$ in the Lie algebra defined by the commutator bracket of  $E({\cal G})$.
   For every element $P\in Spec({\cal G})$,
  we denote by $p_P:{\cal G}\rightarrow {\cal G}/P$ the canonical projection, and by
  $E(p_P):E({\cal H})\rightarrow E({\cal G}/P)$ the induced map on the enveloping algebras.
Let $D$ be a subset of $Spec({\cal G})$, we denote by $\Sigma_D$  the subset
of $ E({\cal G})$ such that for every element $P\in D$, and every element $h\in \Sigma_D$,
$E(p_P)(h)\neq 0$. The elements of $\Sigma_D$ induces on $E({\cal G})$ morphisms of
the right $E({\cal G})$-module $E({\cal G})$ by left multiplications.
 We denote by $E({\cal G})_D$ the localization of $E({\cal G})$ by $\Sigma_D$,
 and by $\alpha_D:E({\cal G})\rightarrow E({\cal G})_D$ the inverting morphism.

\medskip

Let $U$ be an open subset of $Spec({\cal G})$, and $O_{\cal G}(U)$ the set of maps
$f:U\rightarrow \prod_{P\in U}E({\cal G})_P$, such that for every $P\in U$, there exists
an open subset $V$ containing $P$, an element $f_V\in E({\cal G})_V$, such that for every
element $Q\in V, f(Q)$, is the image of $f_V$ by the canonical morphism
$E({\cal G})_V\rightarrow E({\cal G})_P$ resulting from the universal property of $E({\cal G})_P$.

\medskip

{\bf Proposition.}

{\it Let $f:{\cal G}\rightarrow{\cal H}$ be a morphism
 of  ${\cal S}$-algebras,
For every prime ideal $P$ of ${\cal H}$, the inverse image $f^{-1}(P)=f'(P)$ is also prime,

The morphism $f':Spec({\cal H})\rightarrow Spec({\cal G})$ is continuous.
}

\medskip

{\bf Proof.}

Firstly we show that if $P$ is a prime ideal of ${\cal H}$, then $f^{-1}(P)$ is also a prime
ideal of ${\cal G}$. The morphism $f$ induces an injective map
$\bar f:{\cal G}/f^{-1}(P)\rightarrow {\cal H}/P$. It results that ${\cal G}/f^{-1}(P)$
does not have zero divisors, since ${\cal H}/P$ does not have zero divisors.

Let $I$ be an ideal of ${\cal H}$, $f_*^{-1}(V(I))=V(f(I))$. This implies that $f$ is continuous.

\bigskip

{\bf Definition.}

A Lie ${\cal S}$-scheme is a topological space $X$, endowed with a sheaf $O_X$
 such that for every element $x$ of $X$,
there exists a neighborhood $U_x$ of $x$ such that $(U_x,O_X(U_x))$ homeomorphic to an affine
${\cal S}$-scheme $Spec({\cal G}, O_{\cal G})$.

\bigskip

\bigskip

\centerline{\bf References.}

\bigskip

1. Angeleri Hugel L. Archetti M. Tilting modules and universal localizations.
Arxiv.org

\smallskip

 2. Baumslag G,  Myasnikov A, Remeslennikov V, Algebraic geometry over groups I.
 Algebraic sets and ideal theory. Journal of Algebra 219, 16-79 1999.

\smallskip

3. Daniyarova E,  Myasnikov A,  Remeslennikov V.
Algebraic geometry over algebraic structures I, II, III arxiv.org

\smallskip

4. Dieudonn\'e, J. La g\'eom\'etrie des groupes classiques (1963)

\smallskip

5. Hartshorne, R. Algebraic Geometry, Graduate text in Mathematics 52 1977

\smallskip

6. Dieudonn\'e, J Grothendieck, A. \'El\'ements de g\'eom\'etrie alg\'ebrique.

\smallskip

7. Kazachkov V. Algebraic geometry over Lie algebras, Arxiv.org

\smallskip

8. Kharlamapovich O, Myasnikov A. Irreducible affine varieties over free groups I, II
Journal of Algebra 200  517-570, 518-533.

\smallskip

9. Serre, J-P. Arbres, Amalgames, $SL_2$ Ast\'erisque 46 1977

\end{document}